# Structure-Adaptive, Variance-Reduced, and Accelerated Stochastic Optimization

**Junqi Tang**                                                J.TANG@ED.AC.UK
**Francis Bach**                                  FRANCIS.BACH@INRIA.FR
**Mohammad Golbabaee**                            M.GOLBABAEE@ED.AC.UK
**Mike Davies**                                       MIKE.DAVIES@ED.AC.UK
*School of Engineering*
*University of Edinburgh, Edinburgh, UK*
*INRIA - Sierra Project-team*
*Département d'Informatique de l'Ecole Normale Supérieure (CNRS - ENS - INRIA)*
*Paris, France*

## Abstract

In this work we explore the fundamental structure-adaptiveness of state of the art randomized first order algorithms on regularized empirical risk minimization tasks, where the solution has intrinsic low-dimensional structure (such as sparsity and low-rank). Such structure is often enforced by non-smooth regularization or constraints. We start by establishing the fast linear convergence rate of the SAGA algorithm on non-strongly-convex objectives with convex constraints, via an argument of cone-restricted strong convexity. Then for the composite minimization task with a coordinate-wise separable convex regularization term, we propose and analyse a two stage accelerated coordinate descend algorithm (Two-Stage APCG). We provide the convergence analysis showing that the proposed method has a global convergence in general and enjoys a local accelerated linear convergence rate with respect to the low-dimensional structure of the solution. Then based on this convergence result, we proposed an adaptive variant of the two-stage APCG method which does not need to foreknow the restricted strong convexity beforehand, but estimate it on the fly. In numerical experiments we compare the adaptive two-stage APCG with various state of the art variance-reduced stochastic gradient methods on sparse regression tasks, and demonstrate the effectiveness of our approach.

## 1. Introduction

Consider the composite minimization task which reads:

$$x^\star \in \arg\min_{x \in \mathbb{R}^d} \{F(x) := f(x) + \lambda g(x)\}, \tag{1}$$

where we denote $f(x) = \frac{1}{n}\sum_{i=1}^n f_i(x)$ the data fidelity term. Each $f_i(x)$ is convex and $L$-smooth, while the regularization term $g(x)$ is a simple convex function and is possibly non-smooth.

### 1.1 Stochastic variance-reduced optimization and its acceleration

If the objective function $F(x)$ is $\mu$-strongly-convex, stochastic gradient methods with recently introduced variance-reduction techniques SAG (Roux et al., 2012), SDCA (Shamir and Zhang, 2013), SVRG (Johnson and Zhang, 2013), SAGA (Defazio et al., 2014) have a linear convergence rate:

$$\mathcal{O}\left((n + \frac{L}{\mu})\log\frac{1}{\epsilon}\right). \tag{2}$$

Later researchers leveraged *Nesterov acceleration* (Nesterov, 1983)(Nesterov, 2007) with variance-reduced algorithms to achieve an accelerated convergence:

$$\mathcal{O}\left((n + \sqrt{\frac{nL}{\mu}})\log\frac{1}{\epsilon}\right). \tag{3}$$

This line of research starts with the accelerated SDCA (Shalev-Shwartz and Zhang, 2014) which is a dual-based algorithm, and the *accelerated coordinate descent* - APCG (Lin et al., 2014), the primal-dual coordinate descent method SPDC (Zhang and Lin, 2015), RPDG (Lan and Zhou, 2015); then comes the generic acceleration scheme of Catalyst (Lin et al., 2015) which wraps any non-accelerated method into an accelerated proximal point framework to achieve an accelerated rate.

Direct acceleration on stochastic gradients should usually be more tricky than accelerating coordinate descent methods since coordinate descend is guaranteed to decrease the cost function in each iteration. The first attempt is done by (Nitanda, 2014) which directly applies the Nesterov momentum on SVRG which achieves an accelerated rate only when the minibatch-size is large enough; nevertheless later several directly accelerated stochastic gradient methods with variance-reduction have been successfully designed such as Katyusha (Allen-Zhu, 2016) and Point-SAGA (Defazio, 2016).

## 1.2 Restricted strong-convexity, sparsity, and faster convergence

Another line of work (Agarwal et al., 2010) (Agarwal et al., 2012)(Oymak et al., 2015) focuses on the sharpness of the convergence speed of first-order methods with respect to a sharper form of strong-convexity assumption, named *restricted strong convexity* (RSC), towards the optimum $x^\star$, in a restricted set of directions $\mathcal{C}_{x^\star}$ because of the influence of the structure-enforcing regularization or constraints. In general, such types of restricted strong convexity can be described in high-level as the following:

$$f(x) - f(x^\star) - \langle\nabla f(x^\star), x - x^\star\rangle \geq \mu_c\|x - x^\star\|_2^2, \quad \forall x \in \mathcal{C}_{x^\star} \in \mathbb{R}^d, \tag{4}$$

while the global definition of the strong-convexity of $f(x)$ is:

$$f(x) - f(y) - \langle\nabla f(y), x - y\rangle \geq \mu_f\|x - y\|_2^2, \quad \forall x, y \in \mathbb{R}^d. \tag{5}$$

Since $\mu_c \geq \mu_f \geq 0$, such a notion gives us better convergence speed guarantees for first-order optimization, even in many cases when $\mu_f = 0$ such as in Lasso when $n < d$ (Karimi et al., 2016).

The quantification of the RSC constant $\mu_c$ is rather interesting. The first result is given by (Agarwal et al., 2012) in a statistical estimation setting for running a composite gradient method to solve (6) with an additional constraint $\Omega$ to restrict early iterations:

$$x^\star = \arg\min_{x \in \Omega}\left\{F(x) := \frac{1}{n}\sum_{i=1}^n f_i(x) + \lambda g(x)\right\}. \tag{6}$$

With additional assumptions on *decomposible regularizers* $g(x)$ and large enough regularization parameter, they establish a global convergence result for gradient methods with only RSC-type assumption.

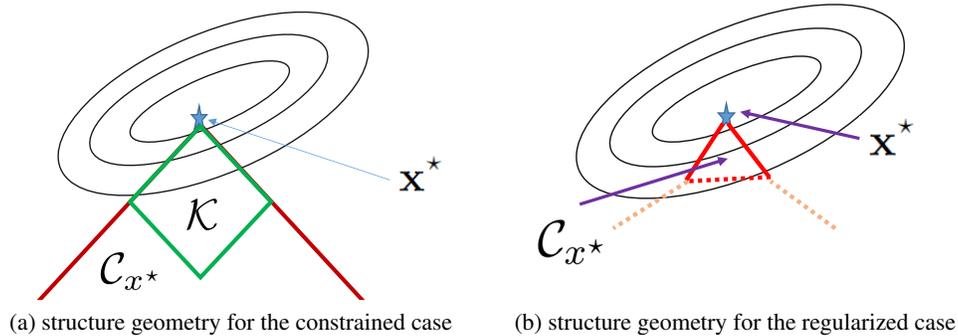

(a) structure geometry for the constrained case  (b) structure geometry for the regularized case

Figure 1: Intuitive view of the geometrical property of two type of empirical risk minimization task with sparsity structure enforced by $l_1$ constraint or regularization: for the $l_1$ constrained case (A), it is straight forward to see that if we run a first order method to find $x^\star$, all the iterates will live in the constrained set $\mathcal{K}$ and hence the descent direction is strictly the cone $\mathcal{C}_{x^\star}$. While for the $l_1$ regularized case (B), as has been shown in the literature, the descent direction is only restrictive nearby the solution $x^\star$

Furthermore, if $f(x) = \frac{1}{2}\|Ax - b\|_2^2$, $A$ is a design matrix which obeys a sub-Gaussian distribution and $g(x) = \lambda\|x\|_1$, the RSC can be quantified w.r.t the sparsity of the true parameter / solution, and the sparser it is, the larger $\mu_c$ will be and hence implies better convergence speed for gradient-base methods. Such results has been recently extended to establish the global linear convergence (Qu and Xu, 2016) (Qu et al., 2017) (Qu and Xu, 2017) of SVRG, SAGA and SDCA on the constrained composite minimization task (6) when only the RSC is available.

An intuitive observation from the details of the theory in (Agarwal et al., 2012) is that for a unconstrained composite minimization task such as Lasso with a correlated Gaussian random design matrix, the RSC is *locally* quantifiable w.r.t sparsity of the solution.

On the other hand, we recall that in the case of the constrained least-squares with an $l_1$ ball as the constrained set:

$$x^\star = \arg\min_{x \in \mathbb{R}^d} \left\{ F(x) := \frac{1}{n}\|Ax - b\|_2^2 + i_c(\|x\|_1 \leq r) \right\}, \tag{7}$$

the RSC is *globally* quantifiable with a Gaussian width statement (Oymak et al., 2015) (Pilanci and Wainwright, 2015) (Pilanci and Wainwright, 2016).

### 1.3 This work

We start by deriving the linear convergence performance of SAGA algorithm (Defazio et al., 2014) on minimizing empirical risk within a convex constraint via cone-restricted strong-convexity. Hence we show that as a paradigm of the variance-reduced stochastic gradient methods, SAGA algorithm automatically adapts to the cone-restricted strong convexity to achieve a linear convergence rate even the lost function is not strongly convex.

In the second part of the paper we go beyond the exact constraint case and turn to the more general setting which is the regularized empirical risk minimization. Meanwhile we go deeper algorithmically and focus on the accelerated methods. We choose to use the accelerated coordinate descent method APCG (Lin et al., 2014) as the foundation to build up our Two-Stage APCG method which is dedicated to actively exploit the intrinsic low-dimension structure of the solution prompted by the (non-smooth) regularization. The convergence analysis shows that the Two-Stage APCG has a global convergence: at the first stage, the method converges sublinearly (linear convergence if we use periodic restart proposed by (Fercoq and Qu, 2016)) to the vicinity of the solution, while in the second stage the method converges towards the solution with an accelerated linear rate with respect to the modified restricted strong convexity (Agarwal et al., 2012) which scales with solution's intrinsic dimension.

In practice the strong convexity and also restricted strong convexity parameter cannot be easily obtained beforehand in general, which is necessary for the accelerated methods to achieve accelerated linear convergence rate. To overcome this we propose an adaptive variant of two-stage APCG method which is based on a simple heuristic scheme to estimate the restricted strong convexity, via a convergence speed check. Our numerical result demonstrates the effectiveness of our algorithm.

## 2. Novel analysis of SAGA algorithm for constrained minimization

### 2.1 SAGA for constrained finite-sum minimization.

We start our work by the analyzing SAGA (Defazio et al., 2014) algorithm's structure-adaptiveness on the constrained minimization task which is a subset of (1):

$$x^\star = \arg\min_{x \in \mathcal{K}} \left\{ f(x) := \frac{1}{n} \sum_{i=1}^{n} f_i(x) \right\}, \tag{8}$$

where the constrained set $\mathcal{K}$ is convex, and we assume that each $f_i(x)$ is convex and has $L$-Lipschitz continuous gradient. The SAGA algorithm is introduced below.

---

**Algorithm 1** SAGA (Defazio et al., 2014) for constrained minimization

**Inputs:** $x^0 \in \mathcal{K}$.
**Initialize:** For each $f_i(.)$, $\phi_i^0 = x^0$ and calculate $\nabla f_i(\phi_i^0)$
**for** $k = 1, \ldots, K$ **do**
    1. Pick an index $j \in [1, n]$ uniformly at random.
    2. Take $\phi_j^{k+1} = x^k$, compute $f_j'(\phi_j^{k+1})$ and store it in the table.
    3. Gradient step using $f_j'(\phi_j^{k+1})$, $f_j'(\phi_j^k)$ and the table average :

$$w^{k+1} = x^k - \gamma \left[ f_j'(\phi_j^{k+1}) - f_j'(\phi_j^k) + \frac{1}{n} \sum_{i=1}^{n} \nabla f_i(\phi_i^k) \right]. \tag{9}$$

    4. Projection step onto convex set $\mathcal{K}$ :

$$x^{k+1} = \mathcal{P}_\mathcal{K}(w^{k+1}) := \arg\min_{u \in \mathcal{K}} \|u - w^{k+1}\|_2^2. \tag{10}$$

**end for**
**Output:** $x^{K+1}$

---

In this setting we can introduce the simplest form of restricted strong convexity assumption on $f(x)$:

**Definition 2.1** *The cone-restricted strong convexity parameter $\mu_c$ of $f(x)$ in (8) is defined as the largest positive constant which satisfies:*

$$f(x) - f(x^\star) - \langle \nabla f(x^\star), x - x^\star \rangle \geq \frac{\mu_c}{2} \|x - x^\star\|_2^2, \quad \forall x \in \mathcal{K}, \tag{11}$$

$$f(x^\star) - f(x) - \langle \nabla f(x), x^\star - x \rangle \geq \frac{\mu_c}{2} \|x - x^\star\|_2^2, \quad \forall x \in \mathcal{K}. \tag{12}$$

An immediate result of (11) is the following lemma (we provide the proof in the appendix):

**Lemma 2.2** *Given the restricted strong convexity in the form of (11) with parameter $\mu_c$, we have:*

$$\langle \nabla f(x) - \nabla f(x^\star), x - x^\star \rangle \geq \mu_c \|x - x^\star\|_2^2, \quad \forall x \in \mathcal{K} \tag{13}$$

*and:*

$$-\|\nabla f(x) - \nabla f(x^\star)\|_2^2 \leq -2\mu_c [f(x) - f(x^\star) - \langle \nabla f(x^\star), x - x^\star \rangle], \quad \forall x \in \mathcal{K} \tag{14}$$

Based on the RSC assumption we are able to derive a new convergence result for SAGA on the constrained finite-sum minimization task.

**Theorem 2.3** *Let $x^\star$ be the optimal solution of (8) and define the Lyapunov function $\mathcal{T}$ as:*

$$\mathcal{T}^k := \mathcal{T}(x^k, \{\phi_i^k\}_{i=1}^n) := \frac{1}{n} \sum_i f_i(\phi_i^k) - f(x^\star) - \frac{1}{n} \sum_i \langle \nabla f_i(x^\star), \phi_i^k - x^\star \rangle + c\|x^k - x^\star\|_2^2 \tag{15}$$

*Then with step size $\gamma = \frac{1}{6L}$, $c = \frac{1}{\gamma(1+4\mu_c\gamma)n}$, the updates of SAGA algorithm obeys:*

$$\mathbb{E}(\mathcal{T}^{k+1}) \leq \left[1 - \min\left(\frac{1}{2n}, \frac{\mu_c}{6L}\right)\right] \mathcal{T}^k \tag{16}$$

We provide the proof of this result in appendix 7.1. Now since each $f_i(.)$ is convex, $\mathcal{T}^k \geq c\|x^k - x^\star\|_2^2$, from Theorem 2.3 we can summarize:

**Corollary 1** *If we run SAGA with step-size $\gamma = \frac{1}{6L}$ to solve (8), the update at k-th iteration satisfies:*

$$\mathbb{E}\|x^k - x^\star\|_2^2 \leq \left[1 - \min\left(\frac{1}{2n}, \frac{\mu_c}{6L}\right)\right]^k \left[\|x^0 - x^\star\|_2^2 + \frac{5n}{18L}[f(x^0) - \langle \nabla f(x^\star), x^0 - x^\star \rangle - f(x^\star)]\right] \tag{17}$$

## 2.2 Minibatch SAGA

From the results given by the previous subsections, we observe (in theory) that, the RSC is enough to guarantee the linear convergence of SAGA, however, when $\mu_c$ is large (e.g. the solution of (8) is very sparse), it does not benefit from the large restricted strong convexity since the linear rate is dominated by $1 - \frac{1}{2n}$. From here we can see that since $n$ is the number of $f_i(.)$, hence in order to exploit the RSC, one many wish to "reduce" $n$, which means using a fixed minibatch scheme in the following form:

$$x^\star = \arg\min_{x \in \mathcal{K}} \left\{ f(x) = \frac{b}{n} \sum_{i=1}^{n/b} \frac{1}{b} \sum_{q=i}^{i+b-1} f_q(x) := \frac{b}{n} \sum_{i=1}^{n/b} f_{b_i}(x) \right\}, \tag{18}$$

where we denote $b$ the minibatch size and assume $\mod(n, b) = 0$ for the simplicity of notation. We now present the (fixed) minibatch SAGA algorithm:

From the same procedure of the previous proof we can have:

**Algorithm 2** Minibatch SAGA for constrained minimization

**Inputs:** $x^0 \in \mathcal{K}$ and minibatch size $b$.
**Initialize:** For each $f_{b_i}(.) := \frac{1}{b} \sum_{q=i}^{i+b-1} f_q(x)$, $\phi_i^0 = x^0$ and calculate $f'_{b_i}(\phi_i^0)$
**for** $k = 1, \ldots, K$ **do**
    1. Pick an index $j \in [1, \frac{n}{b}]$ uniformly at random.
    2. Take $\phi_j^{k+1} = x^k$, compute $f'_{b_j}(\phi_j^{k+1})$ and store it in the table.
    3. Gradient step using $f'_{b_j}(\phi_j^{k+1})$, $f'_{b_j}(\phi_j^k)$ and the table average :

$$w^{k+1} = x^k - \gamma \left[ f'_{b_j}(\phi_j^{k+1}) - f'_{b_j}(\phi_j^k) + \frac{b}{n} \sum_{i=1}^n f'_{b_i}(\phi_i^k) \right]. \tag{19}$$

    4. Projection step onto convex set $\mathcal{K}$ :

$$x^{k+1} = \mathcal{P}_\mathcal{K}(w^{k+1}) := \arg\min_{u \in \mathcal{K}} \|u - w^{k+1}\|_2^2. \tag{20}$$

**end for**
**Output:** $x^{K+1}$

---

**Corollary 2.4** *If each of $f_{b_i}(.)$ has $L_b$-Lipschitz continuous gradient, and we run minibatch SAGA with step-size $\gamma = \frac{1}{6L_b}$ to solve (8), the update at k-th iteration satisfies:*

$$\mathbb{E}\|x^k - x^\star\|_2^2 \leq \left[1 - \min\left(\frac{b}{2n}, \frac{\mu_c}{6L_b}\right)\right]^k \left[\|x^0 - x^\star\|_2^2 + \frac{5n}{18L}[f(x^0) - \langle \nabla f(x^\star), x^0 - x^\star \rangle - f(x^\star)]\right] \tag{21}$$

It is worth noting that the above corollary for minibatch SAGA suggests a linear parallel computation speed up thanks to the cone-restricted strong convexity. Via using minibatches to match $\frac{b}{2n} \approx \frac{\mu_c}{6L_b}$, the overall complexity of SAGA does not have significant change since $(1 - \frac{\mu_c}{6L_b})^{\frac{k}{b}} \approx (1 - \frac{b}{2n})^{\frac{k}{b}} \approx (1 - \frac{1}{2n})^k$ for large $n$ and $k$[1], but this means the larger the cone-restricted strong convexity is, the more parallel speed up minibatch SAGA can achieve. If we choose the minibatch too large such that $\frac{b}{2n} > \frac{\mu_c}{6L_b}$, we lose this linear speed up for parallelisation, and the overall complexity may be worse than SAGA without minibatch.

---

1. Due to the fact that $\lim_{a \to +\infty}(1 - \frac{1}{a})^a = \frac{1}{e}$

## 3. Two-Stage APCG and the convergence analysis

In this section we are going deeper on the analysis of an accelerated coordinate descent algorithm's performance on solving the composite minimization (1) when the strong-convexity only holds in a restricted manner. We first introduce the APCG algorithm (Lin et al., 2014) for (1), for cases where the regularization term $g(x)$ is coordinate-wise separable:

$$g(x) = \sum_{i=1}^{d} g_i([x]_i), \tag{22}$$

and $f(x)$ has coordinate-wise Lipschitz continuous gradient:

$$\|\nabla_i f(x + h_i e_i) - \nabla_i f(x)\|_2 \leq L_i \|h_i\|_2, \forall h_i \in \mathbb{R}, i = 1, ..., d, x \in \mathbb{R}^d. \tag{23}$$

For convenience we define the weighted norm:

$$\|x\|_V = (\sum_{i=1}^{d} L_i \|[x]_i\|_2^2)^{\frac{1}{2}} \tag{24}$$

Direct structure-adaptiveness analysis of APCG could be technically challenging and hence we consider a slightly modified version of it, where we break the iterates into epochs and reset the algorithm (let $z_0^t = x_0^t$) at the beginning of each epoch. More importantly, since the RSC here only holds locally, we introduce the initialization stage given by APCG for non-strongly-convex functions (APCG-NS). (We use superscript $t$ to index outer-loop and subscript $k$ to index inner-loop)

---

**Algorithm 3** Two-Stage APCG

**[Analyzed algorithm]**

**Inputs:** $x_0^0$ and restricted strong-convexity parameter $\mu_c$, number of iteration $K_0$ for the first stage; $T$, $K$ for the outer and inner loop of the second stage.
1. First stage, start without $\mu_c$:

$$x_0^1 = APCG_0(x_0^0, K_0) \tag{25}$$

2. Second stage – exploit local accelerated linear convergence given by $\mu_c$:
**for** $t = 1, \ldots, T$ **do**

$$x_0^{t+1} = APCG(x_0^t, K, \frac{\mu_c}{L}) \tag{26}$$

**end for**
**Output:** $x_0^{T+1}$

**[Implementation]**

**Inputs:** $x_0^0$ and restricted strong-convexity parameter $\mu_c$, number of iterations $K_0$, $N$ for first and second stage.
1. First stage, start without $\mu_c$:

$$x_0^1 = APCG_0(x_0^0, K_0) \tag{27}$$

2. Second stage – exploit local accelerated linear convergence given by $\mu_c$:

$$x_0^2 = APCG(x_0^1, N, \frac{\mu_c}{L}) \tag{28}$$

**Output:** $x_0^2$

---

We list the details of APCG algorithm as the following (Algorithm 4 and 5):

At each iteration, the algorithm chooses a coordinate uniformly at random to perform updates. The update sequences $x_{k+1}^t$ and $z_{k+1}^t$ depend on the realization of the following random variable which we denote as $\xi_k^t$:

$$\xi_k^t = \{i_k^t, i_{k-1}^t, ..., i_1^t, i_0^t, i_k^{t-1}, ..., i_0^{t-1}, ..., i_k^0, ..., i_0^0\}, \tag{35}$$

**Algorithm 4** APCG($x_0^t, K, \mu$) –Accelerated Proximal Coordinate Gradient (Lin et al., 2014, Alg. 2)

**Inputs:** $x_0^t$, number of iteration $K$ and strong-convexity parameter $\mu > 0$.
**Initialize:** $z_0^t = x_0^t, a = \frac{\sqrt{\mu}}{d}$
**for** $k = 1, \ldots, K$ **do**
   1. Compute
$$y_k^t = \frac{x_k^t + a z_k^t}{1+a} \tag{29}$$
   2. Choose $i_k \in 1, ..., d$ uniformly at random and compute
$$z_{k+1}^t = \arg\min_{x \in \mathbb{R}^d} \left[ \frac{da}{2} \|x - (1-a)z_k^t - a y_k^t\|_V^2 + \langle \nabla_{i_k} f(y_k^t), [x]_{i_k} \rangle + \lambda g_{i_k}([x]_{i_k}) \right] \tag{30}$$
   3. Compute
$$x_{k+1}^t = y_k^t + da(z_{k+1}^t - z_k^t) + da^2(z_k^t - y_k^t). \tag{31}$$
**end for**
**Output:** $x_0^{t+1} := x_{K+1}^t$

---

**Algorithm 5** APCG$_0$($x_0^t, K$) –APCG for non-strongly convex functions (Lin et al., 2014, Alg. 3)

**Inputs:** $x_0^t$, number of iteration $K$.
**Initialize:** $z_0^t = x_0^t, a_{-1} = \frac{1}{d}$
**for** $k = 1, \ldots, K$ **do**
   1. Compute
$$a_k = \tfrac{1}{2}(\sqrt{a_{k-1}^4 + 4a_{k-1}^2} - a_{k-1}^2), \quad y_k^t = (1-a_k)x_k^t + a_k z_k^t. \tag{32}$$
   2. Choose $i_k \in 1, ..., d$ uniformly at random and compute
$$z_{k+1}^t = \arg\min_{x \in \mathbb{R}^d} \left[ \frac{da}{2} \|x - z_k^t\|_V^2 + \langle \nabla_{i_k} f(y_k^t), [x]_{i_k} \rangle + \lambda g_{i_k}([x]_{i_k}) \right] \tag{33}$$
   3. Compute
$$x_{k+1}^t = y_k^t + da(z_{k+1}^t - z_k^t). \tag{34}$$
**end for**
**Output:** $x_0^{t+1} := x_{K+1}^t$

---

where $i_k^t$ denotes the index of coordinate the algorithm choose to update at $t$-th outer loop's $k$-th inner loop. For the randomness within a single outer-loop of Two-Stage APCG we specifically denote $\xi_k^t \backslash \xi_k^{t-1}$ as

$$\xi_k^t \backslash \xi_k^{t-1} = \{i_k^t, i_{k-1}^t, ..., i_1^t, i_0^t\} \tag{36}$$

In this section we aim at analyzing the convergence speed of APCG with respect to the low intrinsic dimension of the solution (e.g, sparsity) enforced by the regularization. To achieve this, a special version of RSC (Agarwal et al., 2012) needs to be considered:

$$f(x) - f(x^\star) - \langle \nabla f(x^\star), x - x^\star \rangle \geq \frac{\gamma}{2}\|x - x^\star\|_2^2 - \tau g^2(x - x^\star), \quad \forall x \in \mathbb{R}^d, \tag{37}$$

where $\gamma$ and $\tau$ are some positive constants related to the function $f(.)$ itself. Unlike normal RSC assumptions using Polyak-Lojasiewicz inequality (Karimi et al., 2016), this modified RSC has been shown to have a direct connection with the low-dimensional structure of $x^\star$. As in (Agarwal et al., 2012) we first assume that $g(x)$ is a *decomposable regularizer*:

**Definition 3.1** *(Agarwal et al., 2012) Given a orthogonal subspace pair $(\mathcal{M}, \mathcal{M}^\perp)$ in $\mathbb{R}^d$, $g(.)$ is decomposable if:*

$$g(a+b) = g(a) + g(b), \forall a \in \mathcal{M}, b \in \mathcal{M}^\perp. \tag{38}$$

Meanwhile we define a crucial property for our structure-driven analysis, which is called the *subspace compatibility*:

**Definition 3.2** *(Agarwal et al., 2012) With predefined $g(x)$, we define the subspace compatibility of a model subspace $\mathcal{M}$ as:*

$$\Phi(\mathcal{M}) := \sup_{v \in \mathcal{M} \setminus \{0\}} \frac{g(v)}{\|v\|_2}, \tag{39}$$

*when $\mathcal{M} \neq \{0\}$ and $\Phi(\{0\}) := 0$.*

In this section we also further introduce a structured reference point $x^\dagger$ which is assumed to be nearby $x^\star$. In (Agarwal et al., 2012) this structured point is referred as the ground truth vector or a regression vector in the context of statistical estimation. Throughout our analysis we align the subspace $\mathcal{M}$ such that $x^\dagger \in \mathcal{M}$ or the projection of $x^\dagger$ onto the perturbation subspace $\mathcal{M}^\perp$ is close to zero such that $g(x^\dagger_{\mathcal{M}^\perp})^2$ is small; for the $l_1$ penalized sparse regression, the former case corresponds to the scenario where $x^\dagger$ has exact sparsity while the later corresponds to the approximate sparsity. The subspace compatibility leverages the low-dimensional structure of $x^\star$ into our analysis, for example, if $g(.) = \|.\|_1$, $\|x^\dagger\|_0 = s$ and $\mathcal{M}$ is an $s$-dimensional subspace in $\mathbb{R}^d$, then we have $\Phi(\mathcal{M}) = \sqrt{s}$.

Then we are ready to present the *effective RSC* constant $\mu_c$:

**Lemma 3.3** *(Effective RSC) Given $(x^\star, x^\dagger)$, and denote $\varepsilon := 2\Phi(\mathcal{M})\|x^\dagger - x^\star\|_2 + 4g(x^\dagger_{\mathcal{M}^\perp})$, if the regularization parameter $\lambda$ and the reference point $x^\dagger$ satisfy $\lambda \geq (1 + \frac{1}{c})g^*(\nabla f(x^\dagger))$, then for any convex functions $f(.)$ and $g(.)$ which satisfy:*

$$f(x) - f(x^\star) - \langle \nabla f(x^\star), x - x^\star \rangle \geq \frac{\gamma}{2}\|x - x^\star\|_2^2 - \tau g^2(x - x^\star), \quad \forall x \in \mathbb{R}^d, \tag{40}$$

*with $\gamma > 0$, $\tau > 0$, we have:*

$$f(x) - f(x^\star) - \langle \nabla f(x^\star), x - x^\star \rangle \geq \mu_c \|x - x^\star\|_2^2 - 2\tau(1+c)^2 v^2, \tag{41}$$

*and also:*

$$F(x) - F^\star \geq \mu_c \|x - x^\star\|_2^2 - 2\tau(1+c)^2 v^2, \tag{42}$$

*where $\mu_c = \frac{\gamma}{2} - 8\tau(1+c)^2\Phi^2(\mathcal{M})$ and $v = \frac{\eta}{\lambda} + \varepsilon$, while $x$ satisfies $F(x) - F(x^\star) \leq \eta$.*

---

2. Throughout this paper we denote the Euclidean projection of a vector $v$ onto the subspace $\mathcal{M}$ as $v_\mathcal{M}$ for the simplicity of notation.

**Remark.** The effective RSC $\mu_c = \frac{\gamma}{2} - 8\tau(1+c)^2\Phi^2(\mathcal{M})$ provides us a framework to link the convergence speed of an algorithm with the sparsity of the solution. For example, if $c = 1$, $g(x) = \|x\|_1$ and $\|x^\star\|_0 = s$, then $\Phi^2(\mathcal{M}) = s$ and hence $\mu_c = \frac{\gamma}{2} - 32\tau s$. Further if $F(x)$ is a Lasso problem, then for a wide class of random design matrix we have $\tau = \mathcal{O}(\frac{\log d}{n})$ and $\gamma > 0$. To be more specific, if the data matrix is a correlated Gaussian design matrix such that each row of it is i.i.d drawn from distribution $\mathcal{N}(0, \Sigma)$ where $\Sigma$ is the covariance matrix and we denote its largest and smallest singular value as $r_{\max}(\Sigma)$ and $r_{\min}(\Sigma)$, then it can be shown that $\gamma = \frac{r_{\min}(\Sigma)}{16}$ and $\tau = r_{\max}(\Sigma)\frac{81 \log d}{n}$ with high probability (Raskutti et al., 2010).

The proof of this lemma follows:

**Proof** Let us denote $\Delta = x - x^\dagger$. Since we have assumed $F(x) - F(x^\star) \leq \eta$, then we also have $F(x) - F(x^\dagger) \leq \eta$, hence:

$$f(x^\dagger + \Delta) + \lambda g(x^\dagger + \Delta) \leq f(x^\dagger) + \lambda g(x^\dagger) + \eta, \tag{43}$$

then substract both side with $\langle \nabla f(x^\dagger), \Delta \rangle$ and rearrange:

$$f(x^\dagger + \Delta) - f(x^\dagger) - \langle \nabla f(x^\dagger), \Delta \rangle + \lambda g(x^\dagger + \Delta) - \lambda g(x^\dagger) \leq -\langle \nabla f(x^\dagger), \Delta \rangle + \eta. \tag{44}$$

Due to the convexity of $f(.)$ we immediately have:

$$\begin{aligned}
\lambda g(x^\dagger + \Delta) - \lambda g(x^\dagger) &\leq -\langle \nabla f(x^\dagger), \Delta \rangle + \eta \\
&\leq g^*(\nabla f(x^\dagger))g(\Delta) + \eta \\
&\leq \frac{\lambda}{1+\frac{1}{c}}g(\Delta) + \eta,
\end{aligned}$$

hence by dividing both side with $\lambda$ and then applying the decomposability of $g$ we have:

$$g(x^\dagger + \Delta) - g(x^\dagger) \leq \frac{1}{1+\frac{1}{c}}[g(\Delta_\mathcal{M}) + g(\Delta_{\mathcal{M}^\perp})] + \frac{\eta}{\lambda}, \tag{45}$$

and meanwhile the lower bound on the left-hand-side has been provided in (Agarwal et al., 2012), which reads:

$$g(x^\dagger + \Delta) - g(x^\dagger) \geq g(\Delta_{\mathcal{M}^\perp}) - 2g(x^\dagger_{\mathcal{M}^\perp}) - g(\Delta_\mathcal{M}). \tag{46}$$

By combining these two bounds we have:

$$g(\Delta_{\mathcal{M}^\perp}) + g(\Delta_\mathcal{M}) + \frac{(1+\frac{1}{c})\eta}{\lambda} \leq (1+\frac{1}{c})g(\Delta_{\mathcal{M}^\perp}) - 2(1+\frac{1}{c})g(x^\dagger_{\mathcal{M}^\perp}) - (1+\frac{1}{c})g(\Delta_\mathcal{M}), \tag{47}$$

and then:

$$\begin{aligned}
\frac{1}{c}g(\Delta_{\mathcal{M}^\perp}) &\leq (2+\frac{1}{c})g(\Delta_\mathcal{M}) + 2(1+\frac{1}{c})g(x^\dagger_{\mathcal{M}^\perp}) + \frac{(1+\frac{1}{c})\eta}{\lambda} \\
g(\Delta_{\mathcal{M}^\perp}) &\leq (1+2c)g(\Delta_\mathcal{M}) + 2(1+c)g(x^\dagger_{\mathcal{M}^\perp}) + \frac{(1+c)\eta}{\lambda} \\
g(\Delta) &\leq (2+2c)(g(\Delta_\mathcal{M}) + g(x^\dagger_{\mathcal{M}^\perp})) + \frac{(1+c)\eta}{\lambda}
\end{aligned}$$

Now let $\Delta_x := x - x^\star$ where $x$ satisfies $F(x) - F(x^\star) \leq \eta$, and $\Delta^\star := x^\star - x^\dagger$. Due to the fact that $x^\star$ is the optimal point, $\eta$ can be set as 0 if $x = x^\star$, then:

$$g(\Delta^\star) \leq (2+2c)(g(\Delta^\star_\mathcal{M}) + g(x^\dagger_{\mathcal{M}^\perp})), \tag{48}$$

and now we are able to bound $g(\Delta_x)$:

$$\begin{aligned}
g(\Delta_x) &\leq g(\Delta) + g(\Delta^\star) \\
&\leq (2+2c)g(\Delta_\mathcal{M}) + (2+2c)g(\Delta_\mathcal{M}^\star) + (4+4c)g(x^\dagger_{\mathcal{M}^\perp}) + \frac{(1+c)\eta}{\lambda} \\
&\leq (1+c)\left[2g(\Delta_\mathcal{M}) + 2g(\Delta_\mathcal{M}^\star) + 4g(x^\dagger_{\mathcal{M}^\perp}) + \frac{\eta}{\lambda}\right].
\end{aligned}$$

then by the definition of the subspace compatibility $\Phi(\mathcal{M}) := \sup_{v \in \mathcal{M}\setminus\{0\}} \frac{g(v)}{\|v\|_2}$ we can write:

$$\begin{aligned}
g(\Delta_x) = g(x - x^\star) &\leq (1+c)\left[2\Phi(\mathcal{M})\|x - x^\star\|_2 + 2\Phi(\mathcal{M})\|x^\dagger - x^\star\|_2 + 4g(x^\dagger_{\mathcal{M}^\perp}) + \frac{\eta}{\lambda}\right] \\
&\leq (1+c)\left[2\Phi(\mathcal{M})\|x - x^\star\|_2 + v\right],
\end{aligned}$$

where we denote $\varepsilon := 2\Phi(\mathcal{M})\|x^\dagger - x^\star\|_2 + 4g(x^\dagger_{\mathcal{M}^\perp})$ and $v := \frac{\eta}{\lambda} + \varepsilon$. Then because of the fact that $(a+b)^2 \leq 2a^2 + 2b^2$ we have:

$$g^2(x - x^\star) \leq (1+c)^2 \left[8\Phi^2(\mathcal{M})\|x - x^\star\|_2 + 2v^2\right]. \tag{49}$$

Then we can write:

$$\begin{aligned}
f(x) - f(x^\star) &- \langle \nabla f(x^\star), x - x^\star \rangle \\
&\geq \frac{\gamma}{2}\|x - x^\star\|_2^2 + \tau(1+c)^2\left[8\Phi^2(\mathcal{M})\|x - x^\star\|_2 + 2v^2\right] \\
&\geq \left[\frac{\gamma}{2} - 8\tau(1+c)^2\Phi^2(\mathcal{M})\right]\|x - x^\star\|_2^2 - 2\tau(1+c)^2 v^2,
\end{aligned}$$

which is our first claim. Then because $g(.)$ is convex, we can write:

$$g(x) - g(x^\star) - \langle \partial g(x^\star), x - x^\star \rangle \geq 0, \tag{50}$$

using the fact that the $x^\star$ is the optimal point we justify the second claim. ∎

We assume a non-blowout property of the APCG iterates, which essentially means that the iterates generated by the algorithm will not have a function error too much worse than the error of the starting point:

**Definition 3.4** *(Non-blowout assumption.)* *If we start the APCG algorithm at point $x_0^t$, we assume that there exist a positive constant $1 \leq \omega < \infty$, such that the update sequence $\{x_k^t\}$ generated by the algorithm obeys:*

$$F(x_k^t) - F^\star \leq \omega\left(F(x_0^t) - F^\star\right), \quad \forall t, k \tag{51}$$

We claim that such a relaxed non-blowout assumption is indeed mild and reasonable for APCG. We first recall that the non-accelerated coordinate descent method is guaranteed to not increase the cost function's value in each iteration and hence is strictly non-blowout with $\omega = 1$. Meanwhile (Fercoq and Qu, 2017) has also proved strict non-blowout property with $\omega = 1$ for accelerated full gradient methods. For the moment it is non-trivial to show this property for APCG and hence we temporarily cast it as a relaxed non-blowout assumption. Then we present our key lemma for APCG convergence:

**Lemma 3.5** *Given $(x^\star, x^\dagger)$, and denote $\varepsilon := 2\Phi(\mathcal{M})\|x^\dagger - x^\star\|_2 + 4g(x^\dagger_{\mathcal{M}^\perp})$, if the regularization parameter $\lambda$ and the reference point $x^\dagger$ satisfy $\lambda \geq (1 + \frac{1}{c})g^*(\nabla f(x^\dagger))$. Assume that the non-blowout assumption holds with parameter $\omega$, the updates of the second stage of the Two-Stage APCG obeys:*

$$\mathbb{E}_{\xi_K^t \setminus \xi_K^{t-1}}[F(x_0^{t+1})] - F^\star \leq \left(1 - \frac{\sqrt{\mu_c}}{d\sqrt{L}}\right)^K \cdot 2\left[F(x_0^t) - F^\star\right] + 2\tau(1+c)^2\left(\sqrt{\frac{L}{\mu_c}} + 1\right)v^2, \quad (52)$$

*where $\mu_c = \frac{\gamma}{2} - 8\tau(1+c)^2\Phi^2(\mathcal{M})$, $v = \frac{\eta}{\lambda} + \varepsilon$, $F(x_k^t) - F(x^\star) \leq \eta := \omega\left(F(x_0^t) - F^\star\right)$ for all $t \geq 1$ and $k$, $L = \max_i L_i$.*

**Proof** From the definition of RSC we can have:

$$f(x) - f(x^\star) - \langle \nabla f(x^\star), x - x^\star \rangle$$
$$\geq \mu_c \|x - x^\star\|_2^2 - 2\tau(1+c)^2 v^2$$
$$\geq \frac{\mu_c}{L}\|x - x^\star\|_V^2 - 2\tau(1+c)^2 v^2$$

By observing the main proof of APCG (Lin et al., 2014), we see that there is only one place the strong-convexity assumption on $f(x)$ is used (after equation 3.20). Hence by replacing the original strong-convexity with the effective RSC (41) we have the following:

$$\mathbb{E}_{i_k^t}[f(x_{k+1}^t) + \lambda \hat{g}_{k+1}^t - F^\star + \frac{\mu_c}{2L}\|z_{k+1}^t - x^\star\|_V^2]$$
$$\leq \left(1 - \frac{\sqrt{\mu_c}}{d\sqrt{L}}\right)\mathbb{E}_{i_{k-1}^t}[f(x_k^t) + \lambda \hat{g}_k^t - F^\star + \frac{\mu_c}{2L}\|z_k^t - x^\star\|_V^2] + \frac{2\tau(1+c)^2}{d}v^2,$$

(the detailed definition of $\hat{g}_k^t$ can be found in (Lin et al., 2014, Lemma 3.3), which is a convex combination of $g(z_0^t), g(z_1^t), g(z_2^t) \ldots g(z_k^t)$) and then we roll up the bound:

$$\mathbb{E}_{\xi_k^t \setminus \xi_K^{t-1}}[f(x_{k+1}^t) + \lambda \hat{g}_{k+1}^t - F^\star + \frac{\mu_c}{2L}\|z_{k+1}^t - x^\star\|_V^2]$$
$$\leq \left(1 - \frac{\sqrt{\mu_c}}{d\sqrt{L}}\right)^k [F(x_0^t) - F^\star + \frac{\mu_c}{2}\|x_0^t - x^\star\|_2^2] + \frac{1 - (1 - \sqrt{\mu_c/L}/d)^{k-1}}{1 - (1 - \sqrt{\mu_c/L}/d)} \frac{2\tau(1+c)^2}{d} v^2$$
$$\leq \left(1 - \frac{\sqrt{\mu_c}}{d\sqrt{L}}\right)^k [F(x_0^t) - F^\star + \frac{\mu_c}{2}\|x_0^t - x^\star\|_2^2] + \frac{2\tau(1+c)^2\sqrt{L}}{\sqrt{\mu_c}} v^2$$
$$\leq \left(1 - \frac{\sqrt{\mu_c}}{d\sqrt{L}}\right)^k [2F(x_0^t) - 2F^\star + 2(1+c)^2\tau v^2] + \frac{2\tau(1+c)^2\sqrt{L}}{\sqrt{\mu_c}} v^2$$
$$\leq \left(1 - \frac{\sqrt{\mu_c}}{d\sqrt{L}}\right)^k \cdot 2\left[F(x_0^t) - F^\star\right] + 2\tau(1+c)^2\left(\sqrt{\frac{L}{\mu_c}} + 1\right)v^2$$

where we utilize the effective RSC again to bound the term $\frac{\mu_c}{2}\|x_0^t - x^\star\|_2^2$.

Since $\hat{g}_{k+1}^t \geq g(x_{k+1}^t)$ as declared in (Lin et al., 2014), by simplifying the left hand side we can have:

$$\mathbb{E}_{\xi_k^t \setminus \xi_K^{t-1}}[F(x_{k+1}^t)] - F^\star \leq \left(1 - \frac{\sqrt{\mu_c}}{d\sqrt{L}}\right)^K \cdot 2\left[F(x_0^t) - F^\star\right] + 2\tau(1+c)^2\left(\sqrt{\frac{L}{\mu_c}} + 1\right)v^2. \quad (53)$$

Thus finishes the proof since $F(x_0^{t+1}) = F(x_{K+1}^t)$. ∎

Now we are ready to present our main result:

**Theorem 3.6** *Given $(x^\star, x^\dagger)$, and denote $\varepsilon := 2\Phi(\mathcal{M})\|x^\dagger - x^\star\|_2 + 4g(x^\dagger_{\mathcal{M}^\perp})$, if the regularization parameter $\lambda$ and the reference point $x^\dagger$ satisfy $\lambda \geq (1 + \frac{1}{c})g^*(\nabla f(x^\dagger))$ and we run Two-Stage APCG with $K = \left\lceil \frac{\log 16}{\log \frac{1}{\alpha}} \right\rceil$ and $K_0 = \left\lceil 2\phi d \sqrt{1 + \frac{\|x_0^0 - x^\star\|_V^2}{2[F(x_0^0) - F^\star]}} \right\rceil$, then under the non-blowout assumption with parameter $\omega$, for any $\delta > \lambda\varepsilon$, the update of Two-Stage APCG obeys $F(x_0^t) - F^\star \leq \delta$ if the total number of coordinate-gradient oracle calls $N$ satisfies:*

$$N := tK + K_0 \geq \left\lceil \frac{\log 16}{\log \frac{1}{\alpha}} \right\rceil \log_4 \left( \frac{\frac{1}{\rho\phi^2}[F(x_0^0) - F^\star]}{\delta} \right) + \left\lceil 2\phi d \sqrt{1 + \frac{\|x_0^0 - x^\star\|_V^2}{2[F(x_0^0) - F^\star]}} \right\rceil, \tag{54}$$

*with probability at least $1 - \rho$, where $\alpha := 1 - \frac{\sqrt{\mu_c}}{d\sqrt{L}}$, $\rho = 24\sqrt{\frac{\tau(1+c)^2\omega^2[F(x_0^0)-F^\star]}{\lambda^2\phi^2}}\left(\sqrt{\frac{L}{\mu_c}} + 1\right)$, and $L = \max_i L_i$, $\mu_c = \frac{\gamma}{2} - 8\tau(1+c)^2\Phi^2(\mathcal{M})$.*

**Remarks:**

- Theorem 3.6 presents the main theoretical contribution of this paper. We generalize the structural-analysis framework of (Agarwal et al., 2012) in Lemma 3.3 (that is, we recover their definition of modified RSC by setting $c = 1$) and apply it to our two-stage APCG algorithm. This is the first time in the literature such type of result is shown for randomized first-order method with acceleration.

- The theorem imposes a moderate requirement on the regularization parameter $\lambda$, one is $\lambda \geq (1 + \frac{1}{c})g^*(\nabla f(x^\dagger))$ and the other one is in the probability statement where the $\lambda$ again need to be not too small for meaningful probability. In some sense this should not be a surprise since in a high level argument, if we want to get a meaningful RSC w.r.t the sparsity of the solution the regularization should be large enough to both control the sparsity and also the restricted descent directions.

- The periodic restart here is only needed for the achievable proof and in practice if the $\mu_c$ is known then there is no need to restart the algorithm.

- The parameter $\phi$ denotes the accuracy of the first stage where we run a sub-linearly convergent algorithm. From the proposition we see a very natural consequence – the more iteration the initialization stage takes, the larger $\phi$ is and hence the larger the probability for which the statement is going to hold.

- The contraction factor $\alpha = 1 - \frac{\sqrt{\mu_c}}{d\sqrt{L}}$ occurs in (54) in a logarithmic term $\frac{1}{\log \frac{1}{\alpha}}$ which scales nearly as $\frac{1}{1-\alpha} = \frac{d\sqrt{L}}{\sqrt{\mu_c}}$. Hence we conclude that under the assumptions above, the Two-Stage APCG has a local accelerated linear convergence $\hat{\mathcal{O}}(d\sqrt{\frac{\max_i L_i}{\mu_c}} \log \frac{1}{\delta})$.

**Proof** Following a similar procedure in (Agarwal et al., 2012) (Qu and Xu, 2016) to roll up the residual term $v^2$, we also define three auxiliary sequences $\epsilon_{t+1}/\epsilon_t = \frac{1}{4}$ and $\sigma_{t+1} = 2\sigma_t$, and $v_t = \frac{2\omega\sigma_t\epsilon_t}{\lambda}$ for $t \geq 1$; and $v_0 := F(x_0^0) - F^\star$, $v_1 = \frac{2\omega(F(x_0^1) - F^\star)}{\lambda}$. Since we are focusing on the accuracy regime where $F(x_0^t) - F^\star \geq \lambda\varepsilon$ and the restriction on $\lambda$, such definition of the sequence $v_t$ enable us to apply Lemma 3.5 later, because $\frac{\omega[F(x_0^t)-F^\star]}{\lambda} \geq \varepsilon$. Then according to (Lin et al., 2014) for the first stage of the algorithm we have:

$$\mathbb{E}_{\xi_{K_0}^0}[F(x_{K_0+1}^0)] - F^\star \leq \left(\frac{2d}{2d+K_0}\right)^2 [F(x_0^0) - F^\star + \frac{1}{2}\|x_0^0 - x^\star\|_V^2]$$

Now let $\epsilon_1 = \frac{1}{\phi^2}[F(x_0^0) - F^\star] = \frac{1}{\phi^2}v_0$, then we can choose

$$K_0 := \left\lceil 2\phi d\sqrt{1 + \frac{\|x_0^0 - x^\star\|_V^2}{2[F(x_0^0) - F^\star]}} \right\rceil, \tag{55}$$

we have $\mathbb{E}_{\xi_K^0} F(x_0^1) - F^\star = \mathbb{E}_{\xi_K^0} F(x_{K_0+1}^0) - F^\star \leq \epsilon_1$, and by Markov inequality, with probability at least $1 - \frac{1}{\sigma_1}$ we have:

$$F(x_0^1) - F^\star \leq \sigma_1 \epsilon_1. \tag{56}$$

Next we derive the complexity of the second stage. By an induction statement, we are going to demonstrate that if we choose:

$$\sigma_1 = \sqrt{\frac{\lambda^2 \phi^2}{64\tau(1+c)^2\omega^2(\sqrt{L/\mu_c}+1)v_0}}, \quad K = \left\lceil \frac{\log 16}{\log \frac{1}{\alpha}} \right\rceil, \tag{57}$$

then we will have:

$$\mathbb{E}_{\xi_K^t}(F(x_0^{t+1}) - F^\star) \leq \frac{\epsilon_t}{4}, \quad \forall t \geq 1, \tag{58}$$

where $\epsilon_{t+1}/\epsilon_t = \frac{1}{4}$ and $\epsilon_1 = \frac{1}{\phi^2}[F(x_0^0) - F^\star] = \frac{1}{\phi^2}v_0$, with probability at least $1 - \sum_{i=1}^t \frac{1}{\sigma_i} \geq 1 - \frac{2}{\sigma_1}$:

**Induction part 1:** We turn to our first outer iteration in the second stage of the algorithm. by Lemma 3.5 we have:

$$\mathbb{E}_{\xi_K^1 \setminus \xi_K^0}[F(x_0^2)] - F^\star \leq \alpha^K \cdot 2(F(x_0^1) - F^\star) + 2\tau(1+c)^2\left(\sqrt{\frac{L}{\mu_c}}+1\right)v_1^2, \tag{59}$$

and then we take expectation over $\xi_K^0$:

$$\mathbb{E}_{\xi_K^1}(F(x_0^2) - F^\star) \leq \alpha^K \cdot 2\mathbb{E}_{\xi_K^0}(F(x_0^1) - F^\star) + 2\tau(1+c)^2\left(\sqrt{\frac{L}{\mu_c}}+1\right)v_1^2, \tag{60}$$

where we need:

$$2\tau(1+c)^2\left(\sqrt{\frac{L}{\mu_c}}+1\right)v_1^2 \leq \frac{\epsilon_1}{8}, \tag{61}$$

note that $v_1 = \frac{2\omega(F(x_0^1)-F^\star)}{\lambda} \leq \frac{2\omega\sigma_1\epsilon_1}{\lambda}$ it is enough to satisfy:

$$2\tau(1+c)^2\left(\sqrt{\frac{L}{\mu_c}}+1\right)\left(\frac{2\omega\sigma_1\epsilon_1}{\lambda}\right)^2 \leq \frac{\epsilon_1}{8} \tag{62}$$

which will be satisfied if (recall that $\epsilon_1 = \frac{1}{\phi^2} v_0$):

$$\sigma_1 = \sqrt{\frac{\lambda^2}{64\tau(1+c)^2\omega^2\epsilon_1(\sqrt{L/\mu_c}+1)}} = \sqrt{\frac{\lambda^2\phi^2}{64\tau(1+c)^2\omega^2(\sqrt{L/\mu_c}+1)v_0}} \qquad (63)$$

Then if we choose:

$$K = \left\lceil \frac{\log 16}{\log \frac{1}{\alpha}} \right\rceil, \qquad (64)$$

we can ensure that:

$$\mathbb{E}_{\xi_K^1}(F(x_0^2) - F^\star) \leq \frac{\epsilon_1}{8} + \frac{\epsilon_1}{8} = \frac{\epsilon_1}{4} = \epsilon_2. \qquad (65)$$

**Induction part 2:** For $t+1$-th outer iteration, by induction hypothesis on $t$-th outer iteration which reads: $\mathbb{E}_{\xi_K^{t-1}} F(x_0^t) - F^\star \leq \frac{\epsilon_{t-1}}{4} = \epsilon_t$, we can write:

$$\mathbb{E}_{\xi_K^t \setminus \xi_K^{t-1}}(F(x_0^{t+1}) - F^\star) \leq \alpha^K \cdot 2(F(x_0^t) - F^\star) + 2\tau(1+c)^2 \left(\sqrt{\frac{L}{\mu_c}} + 1\right) v_t^2, \qquad (66)$$

with probability at least $1 - \frac{1}{\sigma_t}$. Then we take expectation over $\xi_K^{t-1}$:

$$\mathbb{E}_{\xi_K^t}(F(x_0^{t+1}) - F^\star) \leq \alpha^K \cdot 2\mathbb{E}_{\xi_K^{t-1}}(F(x_0^t) - F^\star) + 2\tau(1+c)^2 \left(\sqrt{\frac{L}{\mu_c}} + 1\right) v_t^2, \qquad (67)$$

where we need:

$$2\tau(1+c)^2 \left(\sqrt{\frac{L}{\mu_c}} + 1\right) v_t^2 \leq \frac{\epsilon_t}{8}, \qquad (68)$$

since we have chosen that $\sigma_t = 2\sigma_{t-1}$ and $\epsilon_t = \frac{1}{4}\epsilon_{t-1}$, then $v_t = \frac{1}{2}v_{t-1}$, the above inequality is satisfied by our choice of $\sigma_1$.

Again if we choose:

$$K = \left\lceil \frac{\log(16)}{\log \frac{1}{\alpha}} \right\rceil, \qquad (69)$$

we can ensure that:

$$\mathbb{E}_{\xi_K^t}(F(x_0^{t+1}) - F^\star) \leq \frac{\epsilon_t}{8} + \frac{\epsilon_t}{8} = \frac{\epsilon_t}{4} = \epsilon_{t+1}. \qquad (70)$$

Note that:

$$F(x_0^1) - F^\star \leq \frac{F(x_0^0) - F^\star}{\phi^2}, \qquad (71)$$

in summary for Two-Stage APCG if we choose $K := \left\lceil \frac{\log 16}{\log \frac{1}{\alpha}} \right\rceil$, if the number of coordinate gradient oracle calls $N$ satisfies:

$$N := tK + K_0 \geq \left\lceil \frac{\log 16}{\log \frac{1}{\alpha}} \right\rceil \log_4(\frac{\sigma_1[F(x_0^0) - F^\star]}{\phi^2 \delta}) + K_0, \qquad (72)$$

$\mathbb{E}_{\xi_K^{t-1}} F(x_0^t) - F^\star \leq \frac{\delta}{\sigma_1}$ and $F(x_0^t) - F^\star \leq \delta$ with probability at least $1 - \frac{1}{\sigma_1} - \sum_{i=1}^t \frac{1}{\sigma_i} \geq 1 - \frac{3}{\sigma_1}$. Then by setting $\rho = \frac{3}{\sigma_1}$ we finish the proof. ∎

## 3.1 A potentially better result via applying restart on the first stage

The Two-Stage APCG algorithm analyzed in the previous section consists by two stages: in the first stage the algorithm runs in an accelerated sub-linear rate to draw close enough to the optimum, then the second stage algorithm approaches to the optimum in an accelerated linear rate w.r.t the RSC. From the proposition we see that the trade-off between the effort we spend in the first stage and second stage is reflected in the parameter $\phi$, which is related to the approximation accuracy of the first stage.

---

**Algorithm 6** Two-Stage APCG+

**[Analyzed algorithm]**

**Inputs:** $x_0^0$ and restricted strong-convexity parameter $\mu_c$.
**Initialize:** Choose $\mu_0 \in (0,1]$, and set $K_0 = \left\lceil 2\sqrt{3}d\sqrt{1+\frac{1}{\mu_0}} - 2d + 1 \right\rceil$.
1. First stage, start without $\mu_c$:
**for** $j = 1, \ldots, J$ **do**

$$\hat{x}^j = APCG_0(\hat{x}^{j-1}, K_0) \qquad (73)$$

**end for**
$x_0^1 = \hat{x}^J;$
2. Second stage – exploit local accelerated linear convergence given by $\mu_c$:
**for** $t = 1, \ldots, T$ **do**

$$x_0^{t+1} = APCG(x_0^t, K, \tfrac{\mu_c}{L}) \qquad (74)$$

**end for**
**Output:** $x_0^{T+1}$

**[Implementation]**

**Inputs:** $x_0^0$ and restricted strong-convexity parameter $\mu_c$.
**Initialize:** Choose $\mu_0 \in (0,1]$, and set $K_0 = \left\lceil 2\sqrt{3}d\sqrt{1+\frac{1}{\mu_0}} - 2d + 1 \right\rceil$
1. First stage, start without $\mu_c$:
**for** $j = 1, \ldots, J$ **do**

$$\hat{x}^j = APCG_0(\hat{x}^{j-1}, K_0) \qquad (75)$$

**end for**
$x_0^1 = \hat{x}^J;$
2. Second stage – exploit local accelerated linear convergence given by $\mu_c$:

$$x_0^2 = APCG(x_0^1, N, \tfrac{\mu_c}{L}) \qquad (76)$$

**Output:** $x_0^2$

---

In this section we seek to improve the convergence result when the cost function satisfies a weaker form of strong-convexity, by applying restart at the first stage. We introduce the quadratic growth assumption:

**Definition 3.7** *The Quadratic Growth of $F(x)$ is defined as:*

$$F(x) - F^\star \geq \frac{\mu_F}{2}\|x - x^\star\|_2^2, \forall x \in \mathbb{R}^d. \qquad (77)$$

From already established result in the literature, we know that if we restart the $APCG_0$ periodically (if quadratic growth $\mu_F > 0$) we can actually turn it into a linearly convergent algorithm (Fercoq and Qu, 2016). Using this idea we provide an alternative scheme which may lead to a potentially better trade-off on $\phi$, when $\mu_F$ is not too small and $\mu_0$ is a good estimate of $\mu_F$. Note that the first stage of two-stage APCG+ is an instance of (Fercoq and Qu, 2016, Alg. 4).

**Proposition 3.8** *Given $(x^\star, x^\dagger)$, and denote $\varepsilon := 2\Phi(\mathcal{M})\|x^\dagger - x^\star\|_2 + 4g(x_{\mathcal{M}^\perp}^\dagger)$, if the regularization parameter $\lambda$ and the reference point $x^\dagger$ satisfy $\lambda \geq (1 + \frac{1}{c})g^\star(\nabla f(x^\dagger))$ and we run Two-Stage APCG+*

with $K = \left\lceil \frac{\log 16}{\log \frac{1}{\alpha}} \right\rceil$, then under the non-blowout assumption with parameter $\omega$, and also the assumption that for any $x \in \mathbb{R}^d$, $F(x) - F^\star \geq \frac{\mu_F}{2}\|x - x^\star\|_2^2$, for any $\delta > \lambda\varepsilon$, the update of Two-Stage APCG+ obeys $F(x_0^t) - F^\star \leq \delta$ if the total number of coordinate-gradient oracle calls $N$ satisfies:

$$N \geq \left\lceil \frac{\log 16}{\log \frac{1}{\alpha}} \right\rceil \log_4\left(\frac{\frac{1}{\rho\phi^2}\left[F(x_0^0) - F^\star\right]}{\delta}\right)$$

$$+ \left\lceil d\left[6\sqrt{6}\max\left(\frac{1}{\sqrt{\mu_0}}, \frac{L\sqrt{\mu_0}}{\mu_F}\right)\log\left(\frac{\phi^2[F(x_0^0) - F^\star + \frac{d}{2d-2}\|x_0^0 - x^\star\|_V^2]}{F(x_0^0) - F^\star}\right) + 2\sqrt{3}\sqrt{1 + \frac{1}{\mu_0}}\right]\right\rceil,$$

with probability at least $1 - \rho$, where $\alpha := 1 - \frac{\sqrt{\mu_c}}{d\sqrt{L}}$, $\rho = 24\sqrt{\frac{\tau(1+c)^2\omega^2[F(x_0^0) - F^\star]}{\lambda^2\phi^2}}\left(\sqrt{\frac{L}{\mu_c}} + 1\right)$, and $L = \max_i L_i$.

We include the proof of this result in the appendix.

## 4. Adaptive Two-Stage APCG via a simple heuristic procedure for estimating $\mu_c$

In practice, all the state of the art accelerated randomized algorithms for solving the composite minimization task (1) require the explicit knowledge of the strong convexity parameter to run with an accelerated linear convergence rate. For the case where the data fidelity term $f(.)$ is strongly convex, it is difficult in general to calculate the strong convexity parameter before running the accelerated algorithms, let alone in our case, the restricted strong convexity. Inspired by the recent work (Fercoq and Qu, 2017) where an adaptive restart scheme for deterministic full gradient methods based on estimating the quadratic growth (which is a weaker assumption of strong convexity) via checking the convergence speed on the fly with a provable convergence speed guarantee, here we propose an adaptive variant of Two-Stage APCG based on a heuristic procedure for estimating $\mu_c$ on the fly with a small fraction of computational overhead.

Now we describe the intuition of this procedure. First we observe that for $F(x^k) - F^\star < 1$, the convergence speed of the second stage algorithm reads:

$$\mathbb{E}_{\xi_K^t \setminus \xi_K^{t-1}}[F(x_0^{t+1})] - F^\star \leq \left(1 - \frac{\sqrt{\mu_c}}{d\sqrt{L}}\right)^K \cdot 2\left[F(x_0^t) - F^\star\right] + 2\tau(1+c)^2\left(\sqrt{\frac{L}{\mu_c}} + 1\right)v^2$$

$$\approx \left(1 - \frac{\sqrt{\mu_c}}{d\sqrt{L}}\right)^K \cdot 2\left[F(x_0^t) - F^\star\right] + o\left[F(x_0^t) - F^\star\right]$$

$$\approx \left(1 - \frac{\sqrt{\mu_c}}{d\sqrt{L}}\right)^K \cdot 2\left[F(x_0^t) - F^\star\right].$$

Directly using this relationship to check the convergence speed is impossible because $F^\star$ is unknown beforehand in general, but it has been shown in (Fercoq and Qu, 2017, Prop. 4) that $F(x) - F^\star$ can be lower bounded as:

$$F(x) - F^\star \geq \frac{1}{2}\|T(x) - x\|_V^2, \tag{78}$$

**Algorithm 7** Adaptive Two-Stage APCG($x_0, K_0, K, \mu_c, C$)

**Inputs:** $x_0$, number of iteration $K_0$, $K$, an initial guess of the restricted strong-convexity parameter $\mu_c > 0$, and a constant $C \geq 1$.
**Initialize:** $z_0 = x_0$
1. First stage, start without $\mu_c$:
$$x_0 = APCG_0(x_0, K_0) \tag{80}$$

2. Second stage – exploit local accelerated linear convergence given by $\mu_c$:
**for** $k = 1, \ldots, K$ **do**
  1. Set $a_k = \frac{\sqrt{\mu_c}}{d\sqrt{L}}$ and compute
$$y_k = \frac{x_k + a_k z_k}{1+a} \tag{81}$$

  2. Choose $i_k \in 1, \ldots, d$ uniformly at random and compute
$$z_{k+1} = \arg\min_{x \in \mathbb{R}^d} \left[ \frac{da}{2}\|x - (1-a_k)z_k - a_k y_k\|_V^2 + \langle \nabla_{i_k} f(y_k), [x]_{i_k}\rangle + \lambda g_{i_k}([x]_{i_k}) \right] \tag{82}$$

  3. Compute
$$x_{k+1} = y_k + da_k(z_{k+1} - z_k) + da_k^2(z_k - y_k). \tag{83}$$

  4. Modify the estimation of $\mu_c$ every 10 epochs via a convergence speed check:
  **if** mod $(k, 10d) == 0$
    Calculate the composite gradient map:
$$T(x_{k+1}) = \arg\min_{x \in \mathbb{R}^d} \frac{nL}{2}\|x_{k+1} - x\|_2^2 + \langle \nabla f(x_{k+1}), x - x_{k+1}\rangle + \lambda g(x) \tag{84}$$

  **if** $\|T(x_{k+1}) - x_{k+1}\|_2^2 \geq C\left(1 - \frac{\sqrt{\mu_c}}{d\sqrt{L}}\right)^{10d} \|T(x_{k+1-10d}) - x_{k+1-10d}\|_2^2$

$$\mu_c \leftarrow \frac{\mu_c}{2}, \quad a_{k+1} = \frac{\sqrt{\mu_c}}{d\sqrt{L}}, \quad z_{k+1} = x_{k+1}, \tag{85}$$

  **else**
$$\mu_c \leftarrow \min(L, 2\mu_c), \quad a_{k+1} = \frac{\sqrt{\mu_c}}{d\sqrt{L}}, \tag{86}$$

  **end if**
  **if** $a_{k+1} \leq 2^{-5} a_{k+1-50d}$ **then** $C \leftarrow 2C$
  **if** $a_{k+1} \geq 2^5 a_{k+1-50d}$ **then** $C \leftarrow \max(1, \frac{C}{2})$
  **end if**
**end for**
**Output:** $x_{K+1}$

where $T(x)$ is the composite gradient map:
$$T(x) = \arg\min_{q \in \mathbb{R}^d} \frac{nL}{2}\|x - q\|_V^2 + \langle \nabla f(x), q - x\rangle + \lambda g(q), \tag{79}$$

and meanwhile there is upper bound on $F(T(x)) - F^\star$ if we assume there is quadratic-growth $F(x) - F^\star \geq \frac{\mu_F}{2}\|x - x^\star\|_2^2$:

$$F(T(x)) - F^\star \leq \frac{8}{\mu_F}\|T(x) - x\|_V^2, \tag{87}$$

assuming $F(x_0^t) - F^\star = \mathcal{O}[F(T(x_0^t)) - F^\star]$, we may have:

$$\mathbb{E}_{\xi_K^t \setminus \xi_K^{t-1}} \|T(x_0^{t+1}) - x_0^{t+1}\|_V^2 \lesssim \left(1 - \frac{\sqrt{\mu_c}}{d\sqrt{L}}\right)^K \frac{32}{\mu_F} \mathcal{O}[\|T(x_0^t) - x_0^t\|_V^2]$$

Hence our heuristic procedure's checking condition is built based on a simplified version of the above relationship by dropping the expectation:

$$\|T(x_0^{t+1}) - x_0^{t+1}\|_V^2 \lesssim C\left(1 - \frac{\sqrt{\mu_c}}{d\sqrt{L}}\right)^K \|T(x_0^t) - x_0^t\|_V^2 \tag{88}$$

where the variable $C$ represent the strictness of the condition. In the adaptive algorithm we check the condition (88) every a number of epochs, if it is violated we suspect that our estimation of $\mu_c$ is too large and hence we shrink it by a factor of 2 and then restart the second stage algorithm, otherwise we double the estimate to ensure that we choose the estimation of $\mu_c$ as aggressive as possible. If we observe that the algorithm is shrinking the $\mu_c$ for a number of times in a row, we suspect that the algorithm's checking condition is too strict and hence we double $C$ to relax the condition.

## 5. Numerical experiments

This section provides the details of numerical results of our Adaptive Two-Stage APCG algorithm on solving the Lasso regression problem:

$$x^\star \in \arg\min_{x \in \mathbb{R}^d} \left\{ F(x) := \frac{1}{2n}\|Ax - b\|_2^2 + \lambda\|x\|_1 \right\}, \tag{89}$$

We set all our examples with $A \in \mathbb{R}^{n \times d}$ where $n < d$, hence there is no explicit strong convexity in any part of $F(x)$. Since in practice the restricted strong convexity of the data-fidelity term $f(x)$ is usually unknown and difficult to calculate beforehand, in our experiments we apply the adaptive two-stage APCG method which maintains an estimate of $\mu_c$ on the fly and compare it with state of the art variance-reduced stochastic gradient algorithms (proximal-) SVRG (Xiao and Zhang, 2014), and Katyusha (Allen-Zhu, 2016, Algorithm 2) which has an accelerated sub-linear convergence rate for non-strongly convex functions. In particular we also include the vanilla APCG method (Lin et al., 2014, Algorithm 3) as a comparison. All the algorithms in our experiments do not need explicit knowledge of the (restricted) strong convexity of $f(x)$.

For the Million Song data set we choose a subset of samples and manually add random features. This represent the scenario where one may wish to use sparse regression via an $l_1$ penalty to nullify the effect of irrelevant features (Langford et al., 2009). For the Madelon dataset we take the whole training set and add random features as described in the following table. For all the experiments we run our adaptive two-stage APCG algorithm with initial estimate $\mu_c = 0.1$, $C = 1$ and $K_0 = 20d$. For all the algorithms we use step-sizes which have been predicted by the theory.

Table 1: Chosen data sets for Lasso regression

| DATA SET | SIZE $(n, d)$ | # ADDITIONAL RANDOM FEATURES | REFERENCE |
|---|---|---|---|
| YEAR(SUBSET) | (1000, 2070) | 1980 | (LICHMAN, 2013) |
| REGED | (500, 999) | 0 | (WORKBENCH TEAM, 2008B) |
| MADELON | (2000, 4000) | 3500 | (LICHMAN, 2013) |
| MARTI2 | (500, 1024) | 0 | (WORKBENCH TEAM, 2008A) |

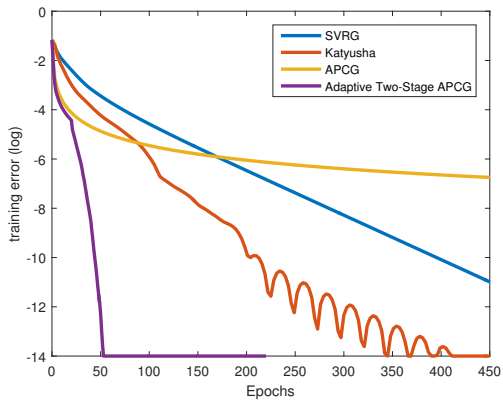

(a) $\lambda = 1 \times 10^{-2}, \|x^\star\|_0 = 39$

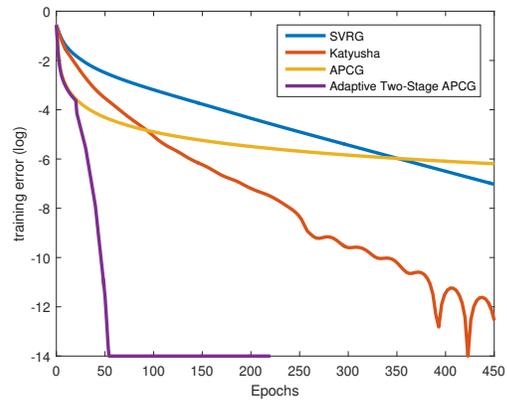

(b) $\lambda = 5 \times 10^{-3}, \|x^\star\|_0 = 292$

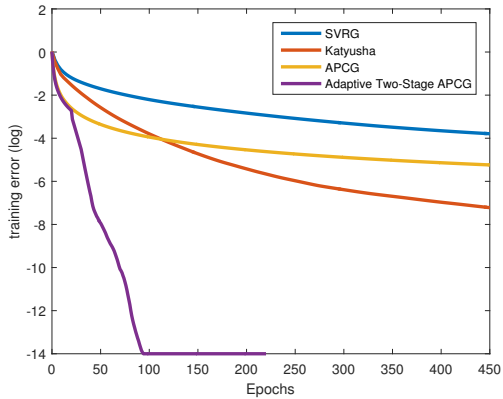

(c) $\lambda = 2 \times 10^{-3}, \|x^\star\|_0 = 636$

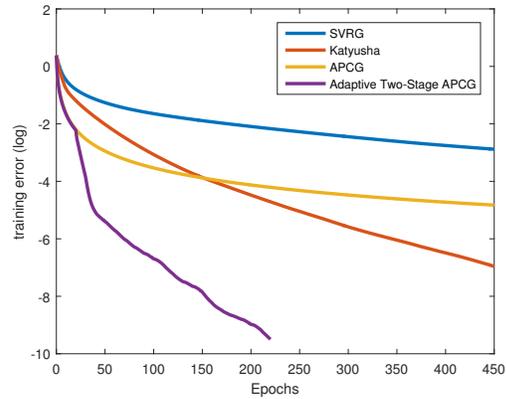

(d) $\lambda = 1 \times 10^{-3}, \|x^\star\|_0 = 821$

Figure 2: Lasso regression on a modified subset of Million-Song Year dataset with additional random features ($A \in \mathbb{R}^{1000 \times 2070}$)

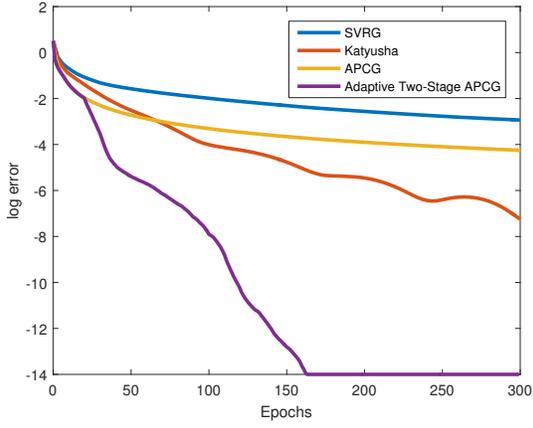
(a) $\lambda = 5 \times 10^{-5}, \|x^\star\|_0 = 34$

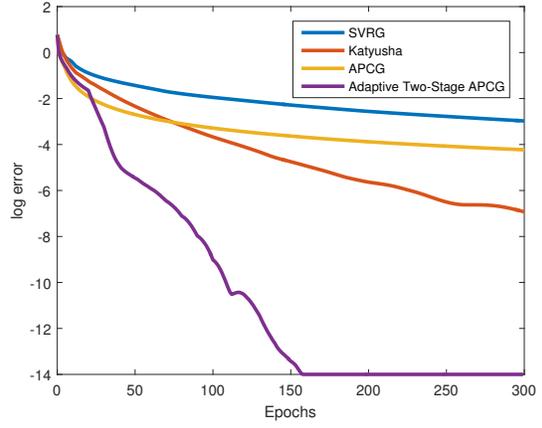
(b) $\lambda = 2 \times 10^{-5}, \|x^\star\|_0 = 80$

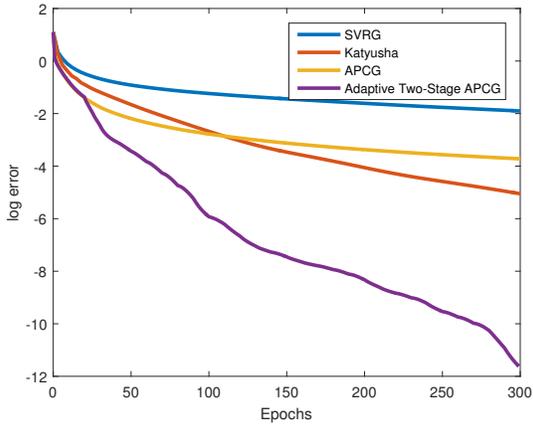
(c) $\lambda = 1 \times 10^{-5}, \|x^\star\|_0 = 127$

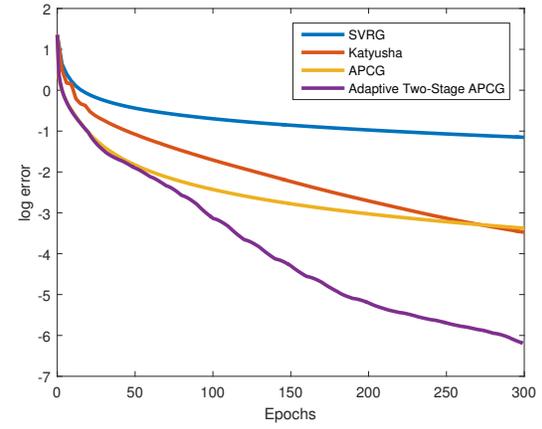
(d) $\lambda = 5 \times 10^{-6}, \|x^\star\|_0 = 209$

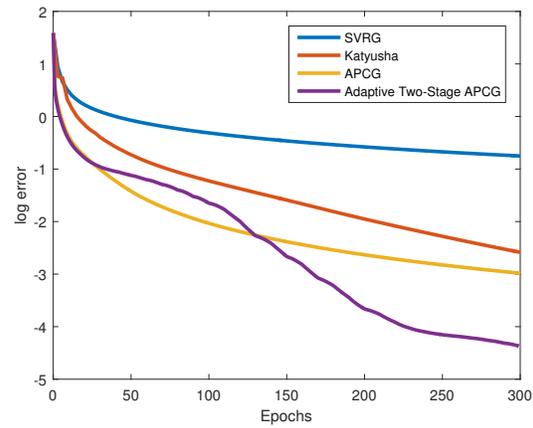
(e) $\lambda = 2 \times 10^{-6}, \|x^\star\|_0 = 343$

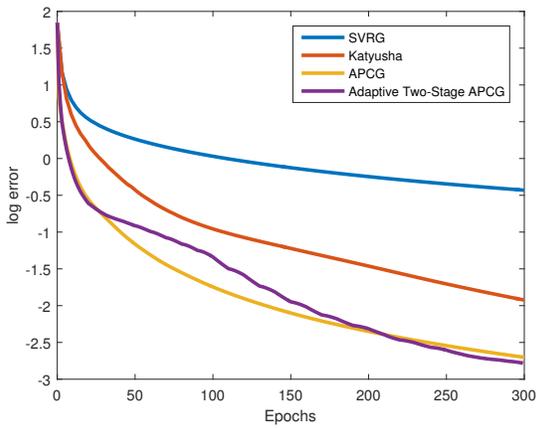
(f) $\lambda = 1 \times 10^{-6}, \|x^\star\|_0 = 395$

Figure 3: Lasso regression on REGED dataset ($A \in \mathbb{R}^{500 \times 999}$)

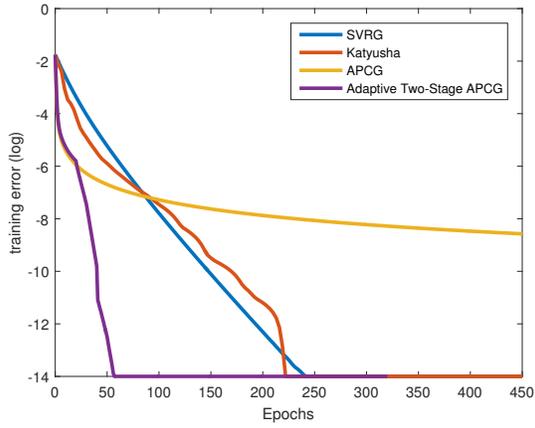
(a) $\lambda = 1 \times 10^{-3}$, $\|x^\star\|_0 = 126$

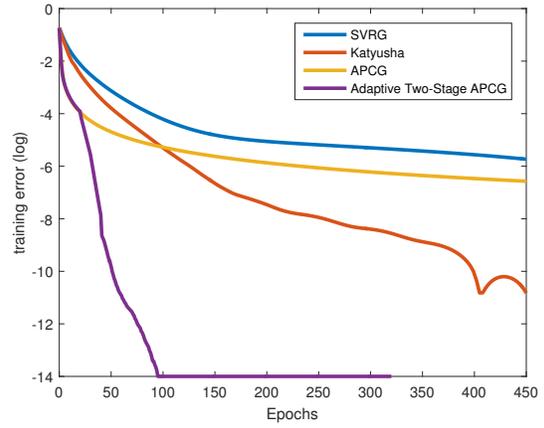
(b) $\lambda = 5 \times 10^{-4}$, $\|x^\star\|_0 = 618$

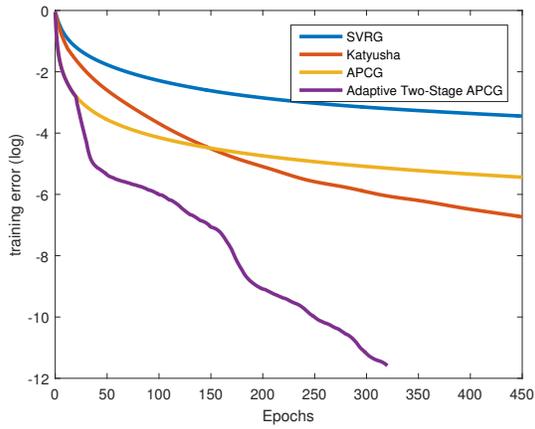
(c) $\lambda = 2 \times 10^{-4}$, $\|x^\star\|_0 = 1250$

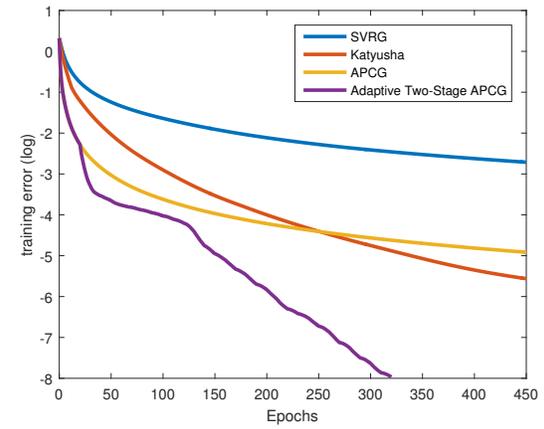
(d) $\lambda = 1 \times 10^{-4}$, $\|x^\star\|_0 = 1594$

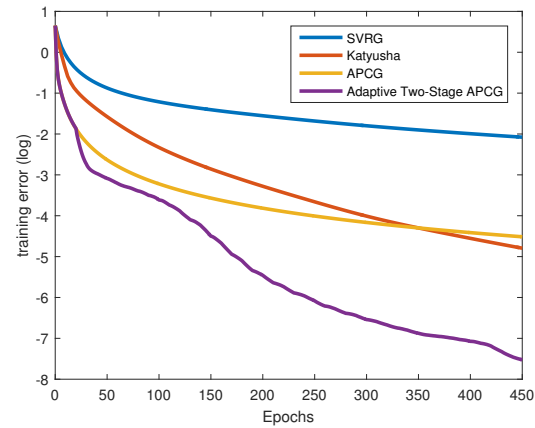
(e) $\lambda = 5 \times 10^{-5}$, $\|x^\star\|_0 = 1777$

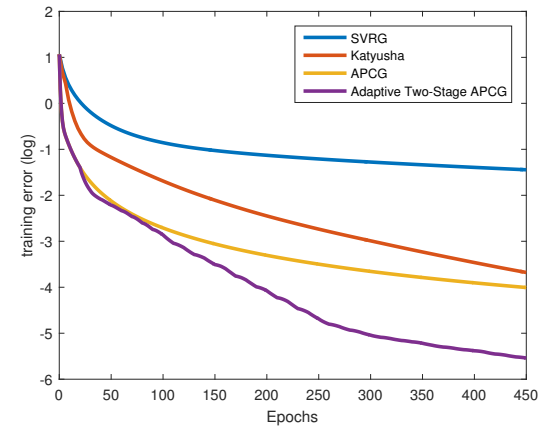
(f) $\lambda = 2 \times 10^{-5}$, $\|x^\star\|_0 = 1912$

Figure 4: Lasso regression on Madelon dataset with additional random features ($A \in \mathbb{R}^{2000 \times 4000}$)

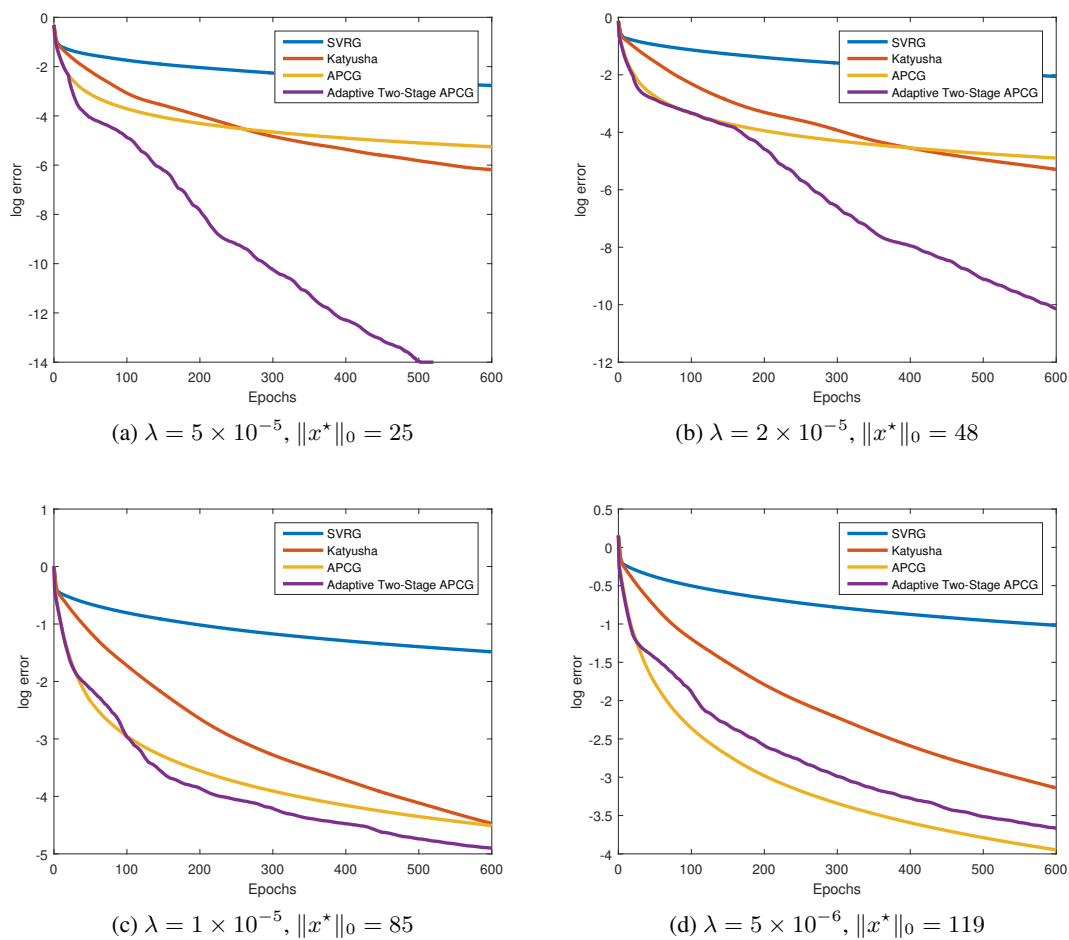

Figure 5: Lasso regression on MARTI2 dataset ($A \in \mathbb{R}^{500 \times 1024}$)

These numerical results on real data sets have demonstrated the effectiveness of our approach for accelerating the APCG method via actively exploiting the low dimensional structure of the solution. The non-adaptive accelerated method like Katyusha and APCG are blind to the potential for the local accelerated linear rate due to the restricted strong convexity. Hence when the solution is relatively sparse, or rather, the regularization parameter is relatively large for the data set, the adaptive two-stage APCG algorithm appears to enjoy local linear convergence speed and significantly outperforms these baselines.

It is worth noting that there is a phase transition phenomenon for our method's performance in all the experiments, when the solution is not sparse enough (or rather, the regularization is not strong enough), this local linear convergence indeed disappears, exactly as predicted by our theory. Such phase transition occurs in various sparsity level which is dependent on the data set itself. For the first three datasets (Year, REGED, Madelon) this phase transition appears only at a trival sparsity level, but for the fourth dataset (Marti2) we only observe accelerated linear convergence for our method when the solution is very sparse.

## 6. Acknowledgements

JT, FB, MG and MD would like to acknowledge the support from H2020-MSCA-ITN Machine Sensing Training Network (MacSeNet), project 642685; ERC grant SEQUOIA; EPSRC Compressed Quantitative MRI grant, number EP/M019802/1; and ERC Advanced grant, project 694888, C-SENSE, respectively. MD is also supported by a Royal Society Wolfson Research Merit Award. JT would like to thank Damien Scieur and Vincent Roulet for helpful discussions during his research visit in SIERRA team.

## 7. Appendix

### 7.1 The proof for Thm 2.3 (SAGA convergence result w.r.t. RSC)

**Proof** This proof relies on our lemma 2.2 and the basic lemmas 1-3 in original SAGA paper (Defazio et al., 2014) and follows the same spirit of the main theorem's proof. We start by boudning each part of $\mathcal{T}^{k+1}$:

$$\mathbb{E}\left[\frac{1}{n}\sum_i f_i(\phi_i^{k+1})\right] = \frac{1}{n}f(x^k) + (1-\frac{1}{n})\frac{1}{n}\sum_i f_i(\phi_i^k); \tag{90}$$

$$\mathbb{E}\left[-\frac{1}{n}\sum_i \langle \nabla f_i(x^\star), \phi_i^{k+1} - x^\star\rangle\right] = -\frac{1}{n}\langle \nabla f(x^\star), x^k - x^\star\rangle - (1-\frac{1}{n})\frac{1}{n}\sum_i \langle \nabla f_i(x^\star), \phi_i^k - x^\star\rangle; \tag{91}$$

and for the third term:

$$\begin{aligned}
&c\mathbb{E}\|x^{k+1} - x^\star\|_2^2 \\
&= c\mathbb{E}\|\mathcal{P}_\mathcal{K}(w^{k+1}) - \mathcal{P}_\mathcal{K}(x^\star - \gamma\nabla\nabla f(x^\star))\|_2^2 \\
&\leq c\mathbb{E}\|w^{k+1} - x^\star + \gamma\nabla f(x^\star)\|_2^2 \\
&= c\mathbb{E}\|x^k - x^\star + w^{k+1} - x^k + \gamma\nabla f(x^\star)\|_2^2 \\
&= c\|x^k - x^\star\|_2^2 - 2c\gamma\langle\nabla f(x^k) - \nabla f(x^\star), x^k - x^\star\rangle + c\mathbb{E}\|w^{k+1} - x^k + \gamma\nabla f(x^\star)\|_2^2
\end{aligned}$$

Up to now the analysis is exactly the same as in (Defazio et al., 2014), but now we split the second term $-2c\gamma\langle\nabla f(x^k) - \nabla f(x^\star), x^k - x^\star\rangle$ into two part by a constant $\theta \in (0,1)$, and bound one of it by the RSC and the other one by the Lemma 1 in (Defazio et al., 2014):

$$\begin{aligned}
&c\|x^k - x^\star\|_2^2 - 2c\gamma\langle\nabla f(x^k) - \nabla f(x^\star), x^k - x^\star\rangle + c\mathbb{E}\|w^{k+1} - x^k + \gamma\nabla f(x^\star)\|_2^2 \\
&= c\|x^k - x^\star\|_2^2 - 2c\gamma(1-\theta)\langle\nabla f(x^k) - \nabla f(x^\star), x^k - x^\star\rangle - 2c\gamma\theta\langle\nabla f(x^k) - \nabla f(x^\star), x^k - x^\star\rangle \\
&\quad + c\mathbb{E}\|w^{k+1} - x^k + \gamma\nabla f(x^\star)\|_2^2 \\
&\leq c\|x^k - x^\star\|_2^2 - 2c\gamma\mu_c(1-\theta)\|x^k - x^\star\|_2^2 - 2c\gamma\theta\langle\nabla f(x^k), x^k - x^\star\rangle + 2c\gamma\theta\langle\nabla f(x^\star), x^k - x^\star\rangle \\
&\quad + c\mathbb{E}\|w^{k+1} - x^k + \gamma\nabla f(x^\star)\|_2^2 \\
&\leq c\|x^k - x^\star\|_2^2 - 2c\gamma\mu_c(1-\theta)\|x^k - x^\star\|_2^2 - 2c\gamma\theta\langle\nabla f(x^k), x^k - x^\star\rangle + 2c\gamma\theta\langle\nabla f(x^\star), x^k - x^\star\rangle \\
&\quad - c\gamma^2\beta\|\nabla f(x^k) - \nabla f(x^\star)\|_2^2 + (1+\beta^{-1})c\gamma^2\mathbb{E}\|f'_j(\phi_j^k) - f'_j(x^\star)\|_2^2 + (1+\beta)c\gamma^2\mathbb{E}\|f'_j(x^k) - f'_j(x^\star)\|_2^2.
\end{aligned}$$

Then we apply (Defazio et al., 2014, Lemma 1) (with $\mu = 0$) to bound $-2c\gamma\theta\langle\nabla f(x^k), x^k - x^\star\rangle$ and (Defazio et al., 2014, Lemma 2) to bound $(1+\beta^{-1})c\gamma^2\mathbb{E}\|f'_j(\phi_j^k) - f'_j(x^\star)\|_2^2$:

$$c\|x^{k+1} - x^\star\|_2^2 \leq \left[c - 2c\gamma\mu_c(1-\theta)\right]\|x^k - x^\star\|_2^2 + \left[(1+\beta)c\gamma^2 - \frac{c\gamma\theta}{L}\right]\mathbb{E}\|f'_j(x^k) - f'_j(x^\star)\|_2^2$$
$$-2c\gamma\theta[f(x^k) - f(x^\star) - \langle\nabla f(x^\star), x^k - x^\star\rangle] - c\gamma^2\beta\|\nabla f(x^k) - \nabla f(x^\star)\|_2^2$$
$$+2(1+\beta^{-1})c\gamma^2 L\left[\frac{1}{n}\sum_i f_i(\phi_i^k) - f(x^\star) - \frac{1}{n}\langle\nabla f_i(x^\star), \phi_i^k - x^\star\rangle\right]$$

Now we bound $\|\nabla f(x^k) - \nabla f(x^\star)\|_2^2$ by Lemma 2.2:

$$c\|x^{k+1} - x^\star\|_2^2 \leq \left[c - 2c\gamma\mu_c(1-\theta)\right]\|x^k - x^\star\|_2^2 + \left[(1+\beta)c\gamma^2 - \frac{c\gamma\theta}{L}\right]\mathbb{E}\|f'_j(x^k) - f'_j(x^\star)\|_2^2$$
$$-2c\gamma\theta[f(x^k) - f(x^\star) - \langle\nabla f(x^\star), x^k - x^\star\rangle]$$
$$-2c\gamma^2\beta\mu_c[f(x^k) - f(x^\star) - \langle\nabla f(x^\star), x^k - x^\star\rangle]$$
$$+2(1+\beta^{-1})c\gamma^2 L\left[\frac{1}{n}\sum_i f_i(\phi_i^k) - f(x^\star) - \frac{1}{n}\langle\nabla f_i(x^\star), \phi_i^k - x^\star\rangle\right].$$

At here we are ready to write:

$$\mathbb{E}(\mathcal{T}^{k+1}) - \mathcal{T}^k \leq -\frac{1}{\kappa}\mathcal{T}^k + \left[\frac{1}{\kappa} - 2\gamma\mu_c(1-\theta)\right]c\|x^k - x^\star\|_2^2$$
$$+ \left[\frac{1}{\kappa} + 2(1+\beta^{-1})c\gamma^2 L - \frac{1}{n}\right]\left[\frac{1}{n}\sum_i f_i(\phi_i^k) - f(x^\star) - \frac{1}{n}\langle\nabla f_i(x^\star), \phi_i^k - x^\star\rangle\right]$$
$$+ \left[\frac{1}{n} - 2c\gamma\theta - 2c\gamma^2\beta\mu_c\right][f(x^k) - f(x^\star) - \langle\nabla f(x^\star), x^k - x^\star\rangle]$$
$$+ \left[(1+\beta)c\gamma^2 - \frac{c\gamma\theta}{L}\right]\mathbb{E}\|f'_j(\phi_j^k) - f'_j(x^\star)\|_2^2.$$

We denote that:

$$c_1 = \frac{1}{n} - 2c\gamma\theta - 2c\gamma^2\beta\mu_c$$
$$c_2 = \frac{1}{\kappa} + 2(1+\beta^{-1})c\gamma^2 L - \frac{1}{n}$$
$$c_3 = \frac{1}{\kappa} - 2\gamma\mu_c(1-\theta)$$
$$c_4 = (1+\beta)\gamma - \frac{\theta}{L}$$

The only remaining task now is to determine the choices of parameters to ensure all these terms are not positive. We first set $c_4 = 0$, with $\gamma = \frac{\theta}{3L}$, then $\beta = 2$. Next we turn to $c_1$:

$$c_1 = \frac{1}{n} - 2c\gamma^2\beta\mu_c - 2c\gamma\theta$$
$$= \frac{1}{n} - 2c\gamma(2\gamma\mu_c + \theta)$$

By setting $c_1 = 0$ we will get $c = \frac{1}{2\gamma(\theta+2\mu_c\gamma)n}$. Then we turn to $c_2$:

$$\begin{aligned}
c_2 &= \frac{1}{\kappa} + 2(1+\beta^{-1})c\gamma^2 L - \frac{1}{n} \\
&= \frac{1}{\kappa} + \frac{3L\gamma}{2(\theta+2\mu_c\gamma)n} - \frac{1}{n} \\
&= \frac{1}{\kappa} + \frac{1}{2(1+\frac{2\mu_c}{3L})n} - \frac{1}{n} \\
&\leq \frac{1}{\kappa} + \frac{1}{2n} - \frac{1}{n} \\
&= \frac{1}{\kappa} - \frac{1}{2n},
\end{aligned}$$

and also $c_3$:

$$c_3 = \frac{1}{\kappa} - \frac{2\mu_c\theta(1-\theta)}{3L}. \tag{92}$$

Now if we choose $\frac{1}{\kappa} = \min\left(\frac{1}{2n}, \frac{2\mu_c\theta(1-\theta)}{3L}\right)$ we can ensure both $c_2$ and $c_3$ are non-positive. Then by setting $\theta = 0.5$ we finishes the proof. ∎

### 7.2 The proof for Lemma 2.2 (the consequence of cone-restricted strong convexity)

**Proof** The first claim is straightforwardly obtained by summing up two copies of (11) with $x$ and $x^\star$ exchanged. To proof the second claim we define an auxiliary function $\psi(.)$ with a fix point $x^\star$:

$$\psi(x) := f(x) - \langle \nabla f(x^\star), x \rangle, \quad x \in \mathcal{K} \tag{93}$$

which has the following property:

$$\psi'(x) = \nabla f(x) - \nabla f(x^\star), \quad \psi'(x^\star) = 0. \tag{94}$$

Then we can write:

$$\begin{aligned}
&\psi(x^\star) - \psi(x) - \langle \psi'(x), x^\star - x \rangle \\
&= f(x^\star) - \langle \nabla f(x^\star), x^\star \rangle - f(x) + \langle \nabla f(x^\star), x \rangle - \langle \nabla f(x) - \nabla f(x^\star), x^\star - x \rangle \\
&= f(x^\star) - f(x) - \langle \nabla f(x), x^\star - x \rangle \\
&\geq \frac{\mu_c}{2}\|x - x^\star\|_2^2.
\end{aligned}$$

Then we have:

$$\begin{aligned}
\psi(x^\star) &\geq \psi(x) + \langle \psi'(x), x^\star - x \rangle + \frac{\mu_c}{2}\|x^\star - x\|_2^2 \\
&\geq \min_v \psi(x) + \langle \psi'(x), v - x \rangle + \frac{\mu_c}{2}\|v - x\|_2^2 \\
&= \psi(x) - \frac{1}{2\mu_c}\|\psi'(x)\|_2^2
\end{aligned}$$

Substituting in the definition of $\psi(.)$ yields the second claim. ∎

## 7.3 The proof for Prop. 3.8 (Two-Stage APCG+)

**Proof** Since the first stage here is equivalent to the restarted APPROX presented in (Fercoq and Qu, 2016, Algorithm 4) with a specific parameter choice, by summarizing the convergence result presented in (Fercoq and Qu, 2016, Corollary 1), to achieve $\mathbb{E}_{\xi^0_{JK_0}} F(x_0^1) - F^\star \leq \epsilon_1$, it is enough to make the total number of coordinate gradient oracle calls in the first stage of Two-Stage APCG+ to be:

$$JK_0 \geq \left\lceil d \left[ 6\sqrt{6} \max\left(\frac{1}{\sqrt{\mu_0}}, \frac{L\sqrt{\mu_0}}{\mu_F}\right) \log\left(\frac{F(x_0^0) - F^\star + \frac{d}{2d-2}\|x_0^0 - x^\star\|_V^2}{\epsilon_1}\right) + 2\sqrt{3}\sqrt{1 + \frac{1}{\mu_0}} \right] \right\rceil \tag{95}$$

Let $\epsilon_1 = \frac{1}{\phi^2}[F(x_0^0) - F^\star]$, we have:

$$JK_0 \geq \left\lceil d \left[ 6\sqrt{6} \max\left(\frac{1}{\sqrt{\mu_0}}, \frac{L\sqrt{\mu_0}}{\mu_F}\right) \log\left(\frac{\phi^2[F(x_0^0) - F^\star + \frac{d}{2d-2}\|x_0^0 - x^\star\|_V^2]}{F(x_0^0) - F^\star}\right) + 2\sqrt{3}\sqrt{1 + \frac{1}{\mu_0}} \right] \right\rceil \tag{96}$$

Then by using the same argument of induction for the second stage in the Two-Stage APCG proof we prove the claim. ∎